\documentclass{article}
\usepackage{amssymb}
\usepackage{amsfonts}
\usepackage{amsmath}
\usepackage{tikz-cd}

\begin{document}

\begin{center}
{\Large About split quaternion algebras over quadratic fields and symbol algebras of degree $n$}

\begin{equation*}
\end{equation*}%

Diana SAVIN

\begin{equation*}
\end{equation*}%
\end{center}

\textbf{Abstract. }{\small In this paper we determine sufficient conditions for a quaternion algebra to split over a quadratic field. In the last section of the paper, we find a class of non-split symbol algebras of degree $n$ (where $n$ is a positive integer, $n\geq 3$) over a $p-$ adic field or over a cyclotomic field.}

\smallskip 


\smallskip


\bigskip

\textbf{1. Introduction}%
\bigskip

Let $K$ be a field with char$K\neq 2.$ Let $K^{\ast }=K\backslash \{0\},$ $a,b$ $\in K^{\ast }.$ The quaternion algebra $H_{K}\left(a,b\right)$ is the $K$-algebra with $K$-basis
$\left\{1; e_{1}; e_{2}; e_{3}\right\}$ satisfying the relations: $e^{2}_{1}=a,$ 
$e^{2}_{2}=b,$ $e_{3}=e_{1}\cdot e_{2}=-e_{2}\cdot e_{1.}$\\
Let $n$ be an arbitrary positive integer, $n\geq 3$ and let $\xi $ be a
primitive $n$-th root of unity. Let $K$ be a field with char$K\neq 2,$ char$K$ does not divide $n$ and $\xi $$%
\in $ $K.$ Let $a,b$ $\in K^{\ast }$ and
let $A$ be the algebra over $K$ generated by elements $x$ and $y$ where%
\begin{equation*}
x^{n}=a,y^{n}=b,yx=\xi xy.
\end{equation*}

This algebra is called a \textit{symbol algebra } and it is denoted by $\left( \frac{a,~b}{%
K,\xi }\right).$ For $n=2,$ we obtain the quaternion algebra. 
Quaternion algebras and symbol algebras are central simple algebras of dimension $n^{2}$ over $K$, non-commutative, but associative algebras (see [Mil; 08]).

In this article we find sufficient conditions for a quaternion algebra to split over a quadratic field. In the paper [Sa; 16] we found a class of division quaternion algebra over the quadratic field $\mathbb{Q}\left(i\right)$ ($i^{2}=-1$), respectively a class of division symbol algebra over the cyclotomic field $\mathbb{Q}\left(\xi\right),$ where $\xi$ is a primitive root of order $q$ (prime) of unity. In the last section of this article we generalize these results for symbol algebras of degree $n\geq 3,$ not necessarily prime.

\bigskip

\textbf{2. Preliminaries}%

\bigskip

We recall some results of the theory of cyclotomic
fields, Kummer fields and $p-$ adic fields, associative algebras, which will be used in our paper. 

\smallskip

Let $n$ be an integer, $n\geq 3$ and let $K$ be a field of characteristic prime to $n$ in which $x^{n}- 1$ splits; and let $\xi$ be a primitive $n$ th root of unity. The following lemma (which can be found in [Ca, Fr; 67]) gives information about certain extension of $K.$ 

\smallskip

\textbf{Lemma 2.1.} \textit{If} $a$ \textit{is a non-zero element of} $K,$
\textit{there is a well-defined normal extension} $K\left(\sqrt[n]{a}\right),$ \textit{the splitting field of} $x^{n}-a.$ \textit{If} $\alpha$ \textit{is a root of} $x^{n}=a,$ \textit{there is a map of the Galois group} $G\left(K\left(\sqrt[n]{a}\right)/K\right)$ \textit{into}
$K^{*}$ \textit{given by} $\sigma\longmapsto \sigma \left(\alpha\right)/\alpha;$ \textit{in particular, if} $a$ \textit{is of order} $n$ \textit{in} $K^{*}/\left(K^{*}\right)^{n},$ \textit{the Galois group is cyclic and can be generated by} $\sigma$ \textit{with} $\sigma \left(\alpha\right)=\xi\alpha.$ \textit{Moreover, the discriminant of} $K\left(\sqrt[n]{a}\right)$
\textit{over} $K$ \textit{divides} $n^{n}\cdot a^{n-1};$ $p$ \textit{is unramified if} $p$ $\nmid$ $na.$

\smallskip

Let $A\neq 0$ be a central simple algebra over a field $K.$ We recall that if $A$ is a
finite-dimensional algebra, then $A$ is a division algebra if and only if $A$
is without zero divisors ($x\neq 0,y\neq 0\Rightarrow xy\neq 0$). $A$ is called \textit{split} by $K$ if $A$ is isomorphic with a matrix algebra over $K.$ If $K\subset L$ is a fields extension, we recall that $A$ is called \textit{split} by $L$  if $A\otimes _{K}L$
is a matrix algebra over $L$. The Brauer group (Br($K$), $\cdot$) of $K$ is Br($K$)$=\left\{\left[A\right]| A\ is\  a\  central\  simple\  K-\ algebra\right\},$ where, two classes of central simple algebras are equals $\left[A\right]=\left[B\right]$ if and only if there are two positive integers $r$ and $s$ such that $A\otimes_{K} M_{r}\left(K\right)\cong B\otimes_{K} M_{s}\left(K\right).$ The group operation in Br($K$) is :  "$\cdot$": Br($K$)$\times$Br($K$)$\longrightarrow$Br($K$), 
$\left[A\right] \cdot\left[B\right]=\left[A\otimes_{K}B\right],$ for $\left(\forall\right)$ $\left[A\right],\left[B\right]$$\in$Br($K$) (see [Mil; 08], [Ko; 00]). A result due Albert-Brauer-Hasse-Noether says that for any number field $K,$ the following sequence is exact:
$$0\longrightarrow Br\left(K\right)\longrightarrow \oplus^{}_{v} Br\left(K_{v}\right)\longrightarrow \mathbb{Q} / \mathbb{Z} \longrightarrow 0$$

\textbf{Remark 2.1.} ([Led; 05]). \textit{Let} $n$ \textit{be a positive integer}, $n\geq 3$ \textit{and let} $\xi $ \textit{be a
primitive} $n$-\textit{th root of unity. Let} $K$ \textit{be a field such that} $\xi$$\in$$K,$ $a,b\in$$K^{*}$. \textit{If} $n$ \textit{is prime, then the symbol algebra} $\left( \frac{a,~b}{%
K,\xi }\right)$ \textit{is either split either a division algebra.}

\smallskip

\textbf{Theorem 2.1.} ([Lin; 12]) (Albert-Brauer-Hasse-Noether). \textit{Let} $H_{F}$ \textit{be a quaternion algebra over a number field} $F$ \textit{and let} $K$ \textit{be a quadratic field extension of} $F.$ \textit{Then there is an embedding of} $K$ \textit{into}
$H_{F}$ \textit{if and only if no prime of} $F$ \textit{which ramifies in} $H_{F}$ \textit{splits in} $K.$

\smallskip

\textbf{Proposition 2.1.} ([Ki, Vo; 10]). \textit{Let} $F$ \textit{be a number field and let} $K$ \textit{be a field containing} $F.$ \textit{Let} $H_{F}$ \textit{be a quaternion algebra over} $F.$ \textit{Let} $H_{K} =H_{F}\otimes _{F} K $ \textit{be a quaternion algebra over
} $K.$ \textit{If} $[K : F] = 2,$ \textit{then} $K$ \textit{splits} $H_{F}$ \textit{if and
only if there exists an} $F$-\textit{embedding} $K\hookrightarrow H_{F}.$

\bigskip

\textbf{3. Quaternion algebras which split over quadratic fields}%

\bigskip

Let $p,q$ be two odd prime integers, $p\neq q.$ If a quaternion algebra $H\left(p,q\right)$
splits over $\mathbb{Q},$ of course it splits over each algebraic number fields. It is known that if $K$ is an algebraic number field such that $[K:\mathbb{Q}]$ is odd and $\alpha,\beta$$\in$$\mathbb{Q}^{*},$ then the quaternion algebra $H_{K}\left(\alpha,\beta\right)$ splits if and only if the quaternion algebra $H_{\mathbb{Q}}\left(\alpha,\beta\right)$ splits (see [Lam; 04]). But, when $[K:\mathbb{Q}]$ is even there are quaternion algebras $H\left(\alpha,\beta\right)$ which does not split over $\mathbb{Q},$ but its split over $K.$ For example,
the quaternion algebra $H\left(11,47\right)$ does not split over $\mathbb{Q},$ but it splits over the quadratic field $\mathbb{Q}\left(i\right)$ (where $i^{2}=-1$).\\
We want to determine sufficient conditions for a quaternion algebra $H\left(p,q\right)$ to split over a quadratic field $K=\mathbb{Q}\left(\sqrt{d}\right).$ Let $\mathcal{O}_{K}$ be the ring of integers of $K.$ Since $p$ and $q$ lie in $\mathbb{Q},$ the problem whether $H_{\mathbb{Q}\left(\sqrt{d}\right)}\left(p,q\right)$ is split reduces to whether
$H_{\mathbb{Q}}\left(p,q\right)$ splits under scalar extension to $\mathbb{Q}\left(\sqrt{d}\right).$\\
It is known that, for each prime positive integer $Br\left(\mathbb{Q}_{p}\right)\cong \mathbb{Q}/\mathbb{Z}$ (the isomorphism is $inv_{p}: Br\left(\mathbb{Q}_{p}\right)\rightarrow \mathbb{Q}/\mathbb{Z}$) and for $p=\infty,$ $Br\left(R\right)\cong \mathbb{Z}/2\mathbb{Z}.$\\
 We obtain sufficient conditions for a quaternion algebra $H\left(p,q\right)$ to split over a quadratic field.

\smallskip

\textbf{Theorem 3.1.} \textit{Let} $d\neq 0,1$ \textit{be a free squares integer,} $d\not\equiv 1$ (\textit{mod} $8$) \textit{and let} $p,q$ \textit{be two prime integers,} $q\geq 3,$ $p\neq q.$ \textit{Let} $\mathcal{O}_{K}$ \textit{be the ring of integers of the quadratic field} $K=\mathbb{Q}\left(\sqrt{d}\right)$ \textit{and} $\Delta_{K}$ \textit{be the discriminant of} $K.$ \textit{Then, we have}:\\
i) \textit{if} $p\geq 3$ \textit{and the Legendre symbols} $\left(\frac{\Delta_{K}}{p}\right)\neq 1,$ $\left(\frac{\Delta_{K}}{q}\right)\neq 1,$ \textit{then, the quaternion algebra} $H_{\mathbb{Q}\left(\sqrt{d}\right)}\left(p,q\right)$ \textit{splits};\\
ii) \textit{if} $p=2$ \textit{and the Legendre symbol} $\left(\frac{\Delta_{K}}{q}\right)\neq 1,$ \textit{then, the quaternion algebra} $H_{\mathbb{Q}\left(\sqrt{d}\right)}\left(2,q\right)$ \textit{splits}.

\smallskip

\textbf{Proof.} i) Applying Albert-Brauer-Hasse-Noether theorem,  we obtain the following description of the Brauer group of $\mathbb{Q}$ and of the Brauer group of the quadratic field $\mathbb{Q}\left(\sqrt{d}\right).$

\[
\begin{tikzcd}
0\arrow{r}{} &  Br\left(\mathbb{Q}\right) \arrow{r}{} \arrow[swap]{d}{} &    \oplus^{}_{p} Br\left(\mathbb{Q}_{p}\right) \cong         \left(\oplus^{}_{p}\mathbb{Q}/\mathbb{Z}\right)\oplus \mathbb{Z}/2\mathbb{Z} \arrow{d}{ \oplus_{p}\varphi_{p}\oplus0}\arrow{r}{} & \mathbb{Q}/\mathbb{Z} \arrow{r}{} & 0 \\
0\arrow{r}{} &  Br\left(\mathbb{Q}\left(\sqrt{d}\right)\right)  \arrow{r}{} & \oplus^{}_{P} Br\left(\mathbb{Q}\left(\sqrt{d}\right)_{P}\right) \cong \left(\oplus^{}_{P}\mathbb{Q}/\mathbb{Z}\right) \arrow{r}{} & \mathbb{Q}/\mathbb{Z}  \arrow{r}{} & 0
\end{tikzcd}
\]
where $\varphi_{p}$ is the multiplication by $2$ when there is single $P$$\in$Spec$\left(\mathcal{O}_{K}\right)$ above the ideal $p\mathbb{Z}$
 i.e. $p\mathbb{Z}$ is inert in $\mathcal{O}_{K}$ or $p\mathbb{Z}$ is ramified in $\mathcal{O}_{K},$ respectively $\varphi_{p}$ is the diagonal map $\mathbb{Q}/\mathbb{Z} \rightarrow \mathbb{Q}/\mathbb{Z} \oplus \mathbb{Q}/\mathbb{Z}$ if there are two primes $P,$ $P^{'}$ of $\mathcal{O}_{K}$ above $p\mathbb{Z}$ i.e. $p\mathbb{Z}$ is totally split in 
$\mathcal{O}_{K}.$ Using this description we determine sufficient conditions for a quaternion algebra $H\left(p,q\right)$ to split over a quadratic field $K=\mathbb{Q}\left(\sqrt{d}\right).$

It is known that $\Delta_{K}=d$ (if $d \equiv 1$ (mod $4$)) or $\Delta_{K}=4d$ (if $d \equiv 2,3$ (mod $4$)). Since $\left(\frac{\Delta_{K}}{p}\right)\neq 1,$ $\left(\frac{\Delta_{K}}{q}\right)\neq 1$ it results $\left(\frac{d}{p}\right)=-1$ or $\left(\frac{d}{p}\right)=0,$ respectively
$\left(\frac{d}{q}\right)=-1$ or $\left(\frac{d}{q}\right)=0.$
Applying the theorem of decomposition of a prime integer $p$ in the ring of integers of a quadratic field (see for example [Ire, Ros; 90], p. 190), it results that $p$ is ramified in $\mathcal{O}_{K}$ or $p$ is inert in $\mathcal{O}_{K},$ respectively $q$ is ramified in $\mathcal{O}_{K}$ or $q$ is inert in $\mathcal{O}_{K}.$ So, $p$ and $q$ do not split in $K.$\\
Let $\Delta$ denote the discriminant of the quaternion algebra $H_{\mathbb{Q}\left(\sqrt{d}\right)}\left(p,q\right).$\\
It is known that a prime positive integer $p^{'}$ ramifies in $H_{\mathbb{Q}\left(\sqrt{d}\right)}\left(p,q\right)$ if $p^{'}| 2\Delta$ ([Ko], [Ko; 00]). This implies $p^{'}| 2pq.$\\
Since $d\not\equiv 1$ (mod $8$) and the decomposition of $2$ in $\mathcal{O}_{K}$ (see [Ire, Ros; 90], p. 190), it results that $2$ does not split in $K.$\\
From the previously proved and applying Theorem 2.1 and Proposition 2.1, it results that
the quaternion algebra $H_{\mathbb{Q}\left(\sqrt{d}\right)}\left(p,q\right)$ splits.\\
ii) Let $p^{'}$ be a prime positive integer which ramifies in $H_{\mathbb{Q}\left(\sqrt{d}\right)}\left(2,q\right).$ In this case the condition $p^{'}| 2\Delta$ implies $p^{'}| 2q.$
With similar reasoning as i) we get that the quaternion algebra $H_{\mathbb{Q}\left(\sqrt{d}\right)}\left(2,q\right)$ splits.

\smallskip

\textbf{Remark 3.1.} The conditions $\left(\frac{\Delta_{K}}{p}\right)\neq 1,$ $\left(\frac{\Delta_{K}}{q}\right)\neq 1$ from Theorem 3.1 are not necessary for the quaternion algebra $H_{\mathbb{Q}\left(\sqrt{d}\right)}\left(q,p\right)$ splits. For example, if $d=-1,$ the conditions $\left(\frac{\Delta_{K}}{p}\right)\neq 1,$ $\left(\frac{\Delta_{K}}{q}\right)\neq 1$ are equivalent to $p$$\equiv$$q$$\equiv$$3$(mod $4$). We consider the quaternion algebra $H_{\mathbb{Q}\left(i\right)}\left(5, 29\right),$ so $p=5$$\equiv$$1$(mod $4$) and $q=29$$\equiv$$1$(mod $4$). Doing some calculations in software MAGMA, we obtain that the algebra $H_{\mathbb{Q}\left(i\right)}\left(5, 29\right)$ splits. Analogous, for $p=5$$\equiv$$1$(mod $4$) and $q=19$$\equiv$$3$(mod $4$), we obtain that the algebra $H_{\mathbb{Q}\left(i\right)}\left(5, 19\right)$ splits. Another example: if $d=3,$ $p=7, q=47.$ We have $\left(\frac{\Delta_{K}}{p}\right)\neq 1,$ but $\left(\frac{\Delta_{K}}{q}\right)=1.$ However the quaternion algebra $H_{\mathbb{Q}\left(\sqrt{3}\right)}\left(7, 47\right)$ splits. Another remark is that the quaternion algebra $H_{\mathbb{Q}}\left(7, 47\right)$ does not split.

\smallskip

We wonder what happens with a quaternion algebra $H_{\mathbb{Q}\left(\sqrt{d}\right)}\left(p,q\right)$ from Theorem 3.1 when instead of $p$ or $q$ we consider an arbitrary integer $\alpha.$ Immediately we obtain the following result:

\smallskip

\textbf{Corollary 3.1.} \textit{Let} $d\neq 0,1$ \textit{be a free squares integer,} $d\not\equiv 1$ (\textit{mod} $8$) \textit{and let} $\alpha$ \textit{be an integer and} $p$ \textit{be an odd prime integer.} \textit{Let} $\mathcal{O}_{K}$ \textit{be the ring of integers of the quadratic field} $K=\mathbb{Q}\left(\sqrt{d}\right)$ \textit{and} $\Delta_{K}$ \textit{be the discriminant of} $K.$ \textit{If} \textit{the Legendre symbols} $\left(\frac{\Delta_{K}}{p}\right)\neq 1,$ $\left(\frac{\Delta_{K}}{q}\right)\neq 1,$ \textit{for each odd prime divisor} $q$ \textit{of} $\alpha$ \textit{then, the quaternion algebra} $H_{\mathbb{Q}\left(\sqrt{d}\right)}\left(\alpha,p\right)$ \textit{splits}.

\smallskip

\textbf{Proof.} We want to determine the primes $p^{'}$ which ramifies in $H_{\mathbb{Q}\left(\sqrt{d}\right)}\left(\alpha,p\right),$ i.e the primes $p^{'}$ with the property $p^{'}| 2\Delta.$ This implies $p^{'}| 2 \alpha \cdot p.$ Since $\left(\frac{\Delta_{K}}{p}\right)\neq 1,$ $\left(\frac{\Delta_{K}}{q}\right)\neq 1,$ for each odd prime divisor $q$ \textit{of} $\alpha,$ using a reasoning similar with that of Theorem 3.1, we get that such primes does not exist, so the quaternion algebra $H_{\mathbb{Q}\left(\sqrt{d}\right)}\left(\alpha,p\right)$ splits.\\
\bigskip

\textbf{4. Symbol algebras of degree $n$}%
\bigskip

In the paper [Sa; 16] we found a class of division quaternion algebras over the quadratic field $\mathbb{Q}\left(i\right)$ ([Sa; 16], Th. 3.1) and a class of division symbol algebras of degree $q$ (where $q$ is an odd prime positive integer) over a $p$- adic field or over a cyclotomic field ([Sa; 16], Th. 3.2). Here we generalize theorem 3.2 from [Sa; 16], when $A$ is a symbol algebra over the $n
$-th cyclotomic field, where $n$ is a positive integer, $n\geq 3.$

\smallskip

\textbf{Theorem 4.1.} \textit{Let} $n$ \textit{be a positive integer}, $n\geq 3,$ $p$ \textit{be} \textit{%
a prime positive integer such that} $p\equiv 1$ (\textit{mod} $n$), $\xi $ 
\textit{be a primitive root of order} $n$ \textit{of unity and let }$K=%
\mathbb{Q}\left( \xi \right) $\textit{\ be the} $n$ \textit{th cyclotomic field.} \textit{%
Then there is an integer} $\alpha $ \textit{not divisible by} $p,$ $\alpha $  \textit{is not a}
$l$ \textit{power residue modulo} $p,$ \textit{for} $\left(\forall\right)$ $l$$\in$$\mathbb{N},$ $l|n$ \textit{and for every such an} $\alpha $%
, \textit{we have:}\newline
i) \textit{if} $A$ \textit{is the symbol
algebra} $A=\left( \frac{\alpha ,p}{K,\xi }\right),$ \textit{then} $A\otimes _{K}\mathbb{Q}_{p}$ \textit{is a non-split
algebra over} $\mathbb{Q}_{p};$\newline
ii) \textit{the symbol algebra} $A$ \textit{is a non-split algebra over} $%
K.\smallskip $

\textbf{Proof.} i) Let be the homomorphism $f:\mathbb{F}_{p}^{\ast }\rightarrow 
\mathbb{F}_{p}^{\ast },$ $f\left( x\right) =x^{n}.$ Since $p\equiv 1$ (\textit{mod} $n$),
it results $Ker\left( f\right) =\left\{ x\in {F}_{p}^{\ast }|x^{n}\equiv 1\
(mod\ p)\right\} $ is non -trivial, so $f$ is not injective. So, $f$ is not
surjective. It results that there exists $\overline{\alpha }$ (in $\mathbb{F}_{p}^{\ast
}, $) which does not belongs to $\left( \mathbb{F}_{p}^{\ast }\right) ^{n}.$ Let $\beta$ be an $n$ th root of $\alpha$ (modulo $p$). Since $\alpha $ is not a
$l$ power residue modulo $p,$ for $\left(\forall\right)$ $l$$\in$$\mathbb{N},$ $l|n,$ it results that the extension of fields $\mathbb{F}%
_{p}\left(\overline{\beta}\right) /\mathbb{F}_{p}$ is a cyclic
extension of degree $n.$ Applying a consequence of Hensel's lemma (see for example [Al,Go;99])
and the fact that $p\equiv 1$ (mod $n$), it results that $\mathbb{Q}_{p}$ contains
the $n$-th roots of the unity, therefore $\mathbb{Q}\left( \xi \right)
\subset \mathbb{Q}_{p}.$ Let the symbol algebra $A\otimes _{K}%
\mathbb{Q}_{p}=\left( \frac{\alpha ,p}{\mathbb{Q}_{p},\xi }\right) .$  Applying Lemma 2.1, it results that the
extension $\mathbb{Q}_{p}\left( \sqrt[n]{\alpha }\right) /\mathbb{Q}_{p}$ is
a cyclic unramified extension of degree $n,$ therefore a norm of an element
from this extension can be a positive power of $p,$ but cannot be $p.$
According to a criterion for splitting of  the symbol algebras (see Corollary 4.7.7, p. 112 from [Gi, Sz; 06]), it results that $\left( \frac{\alpha ,p}{\mathbb{Q}%
_{p},\xi }\right) $ is a non-split algebra.\\
ii) Applying i) and the fact that $K\subset \mathbb{Q}_{p},$ it results that $A$ is a non-split algebra.

\smallskip

\textbf{Remark 4.1.} Although Theorem 4.1 is the generalization of Theorem 3.2 from [Sa; 16] for 
symbol algebras of degree $n,$ there are some differences between these two theorems, namely:\\
- One of the conditions of the hypothesis of Theorem 3.2 from [Sa; 16] is: $\alpha $ is not a
$q$ power residue modulo $p.$ With a similar condition in the hypothesis of Theorem 4.1, namely: $\alpha $  is not a $n$ power residue modulo $p,$ Theorem 4.1 does not work. We give an example in this regard: let $p=7,$ $n=6,$ $\alpha=2.$ $2$ is not a $6$ power residue modulo $7,$ but $2$ is a quadratic residue modulo $7.$ Let $\beta$ be an $6$ th root of $\alpha$ (modulo $7$). We obtain that the polynomial $Y^{6}-\overline{2}$ is not irreducible in $\mathbb{F}_{7}\left[Y\right].$ We have $Y^{6}-\overline{2}=\left(Y^{3}-\overline{3}\right)\cdot \left(Y^{3}+\overline{3}\right)$ (in $\mathbb{F}_{7}\left[Y\right]$).
So, the extension of fields $\mathbb{F}_{7}\subset\mathbb{F}%
_{7}\left( \overline{\beta}\right)$ has not the degree $n=6.$ For this reason, in the hypothesis of Theorem 4.1 we put the condition: $\alpha $ is not a
$l$ power residue modulo $p,$ for $\left(\forall\right)$ $l$$\in$$\mathbb{N},$ $l|n;$\\
- In Theorem 3.2 from [Sa; 16] we proved that $A\otimes _{K}\mathbb{Q}_{p}$ is a non-split
symbol algebra over $\mathbb{Q}_{p}$ (respectively $A$ is a non-split
symbol algebra over $K$) and applying Remark 2.1. this is equivalent to $A$ is a division
symbol algebra over $\mathbb{Q}_{p}$ (respectively $A$ is a division symbol algebra over $K$).
But, Remark 2.1 holds iff $n$ is prime.  For this reason, the conclusion of Theorem 4.1 is: $A$ is a non-split symbol algebra over $\mathbb{Q}_{p}$ (respectively $A$ is a non-split symbol algebra over $K$).

\smallskip

\textbf{Conclusions.} In the last section of the paper, we found a class of non-split symbol algebras of degree $n$ (where $n$ is a positive integer, $n\geq 3$) over a $p-$ adic field, respectively over a cyclotomic field. In a further research we intend to improve Theorem 4.1 from this paper, for to find a class of division symbol algebras of degree $n$ (where $n$$\in$$\mathbb{N}^{*}$, $n\geq 3$) over a cyclotomic field.

\bigskip

\textbf{References}%
\begin{equation*}
\end{equation*}
[Al, Go; 99] V. Alexandru, N.M. Gosoniu, \textit{Elements of Number Theory } (in Romanian), Ed.
 Bucharest University, 1999.\newline
[Ca, Fr; 67]  J. W. S. Cassels,  A. Fr$\ddot{o}$hlich (editors), \textit{Algebraic Number Theory (Proceedings of an instructional conference organized by the London Mathematical Society),} Academic Press, 1967. \newline
[Gi, Sz; 06] P. Gille,  T. Szamuely, \textit{Central Simple Algebras and Galois
Cohomology}, Cambridge University Press, 2006.\newline
[Ire, Ros; 90] K. Ireland, M. Rosen, \textit{A classical introduction to modern number theory},
Springer-Verlag, 1990.\newline
[Ki, Vo; 10] M. Kirschmer, J. Voight, \textit{Algorithmic enumeration of ideal classes for quaternion orders,} SIAM J. Comput. (SICOMP) \textbf{39} (2010), no. 5, 1714-1747.\newline
[Ko] D. Kohel,  \textit{Quaternion algebras,} echidna.maths.usyd.edu.au/kohel/alg/doc/\newline
AlgQuat.pdf\newline
[Ko; 00] D. Kohel,  \textit{Hecke module structure of quaternions,} Proceedings of Class Field Theory - Centenary and Prospect (Tokyo, 1998), K. Miyake, ed., Advanced Studies in Pure Mathematics, \textbf{30}, 177-196, 2000. \newline
[Lam; 04] T. Y. Lam, \textit{Introduction to Quadratic Forms over Fields,}
American Mathematical Society, 2004.\newline
[Led; 05] A. Ledet, \textit{Brauer Type Embedding Problems }, American
Mathematical Society, 2005.\newline
[Lin; 12] B. Linowitz, \textit{Selectivity in quaternion algebras}, Journal of  Number Theory \textbf{132} (2012), pp. 1425-1437. \newline   
[Mil; 08] J.S. Milne, \textit{Class Field Theory},
http://www.math.lsa.umich.edu/~jmilne.\newline
[Sa; 16] D. Savin, \textit{About division quaternion algebras and division symbol algebras},
Carpathian Journal of  Mathematics, \textbf{32(2)} (2016), p. 233-240.\newline
[Vo; 10] J. Voight, The Arithmetic of Quaternion Algebras. Available on the
author's website: http://www.math.dartmouth.edu/~jvoight/
crmquat/book/quat-modforms-041310.pdf, 2010.
\[\]

{\footnotesize \vspace{2mm} \noindent 
\begin{minipage}[b]{10cm}
Diana SAVIN, \\
Faculty of Mathematics and Computer Science, Ovidius University, \\ 
Constanta 900527, Bd. Mamaia no.124, Rom\^{a}nia \\
Email: savin.diana@univ-ovidius.ro;  dianet72@yahoo.com
\end{minipage}
}

{\footnotesize \bigskip }

\end{document}